\documentclass[10pt]{amsart}
\usepackage{amsmath,amssymb,amsthm} 
\usepackage{xcolor}
\usepackage[a4paper, total={6in, 8in}]{geometry}
\newtheorem{theorem}{\bf Theorem}[section]
\newtheorem{definition}[theorem]{\bf Definition}
\newtheorem{lemma}[theorem]{\bf Lemma}

\newcommand{\R}{\mathbb{R}}
\newcommand{\Z}{\mathbb{Z}}

\title{\hspace{-12pt}Constructions for the Elekes--Szab\'{o} and Elekes--R\'{o}nyai problems}
\author{Mehdi Makhul, Oliver Roche-Newton, Audie Warren and Frank de Zeeuw}


\begin{document}

\maketitle

\begin{abstract}
We give a construction of a non-degenerate polynomial $F\in \R[x,y,z]$ and a set $A$ of cardinality $n$ such that $\left|Z(F)\cap (A \times A \times A) \right| \gg n^{\frac{3}{2}}$, thus providing a new lower bound construction for the Elekes--Szab\'{o} problem.
We also give a related construction for the Elekes--R\'{o}nyai problem restricted to a subgraph.
This consists of a polynomial $f\in \R[x,y]$ that is not additive or multiplicative, a set $A$ of size $n$, 
and a subset $P\subset A\times A$ of size $|P|\gg n^{3/2}$ on which $f$ takes only $n$ distinct values.
\end{abstract}


\section{Introduction}
Throughout this paper, we write $X \gg Y$ if and only if there exists some absolute constant $c > 0$ such that $X \ge cY$. If the constant $c$ depends on another parameter $k$, we use the shorthand $X \gg_k Y$.


\subsection{The Elekes--Szab\'{o} Problem}

Elekes and Szab\'o \cite{ES12} considered the size of the intersection of the zero set of a polynomial $F(x,y,z)\in \R[x,y,z]$ of degree $d$ with a Cartesian product $A \times B \times C \subset \mathbb{R}^3$, 
where $|A|=|B|=|C|=n$. 
By the Schwartz--Zippel Lemma (see for instance \cite[Lemma A.4]{RSZ16}),
we have
\begin{equation}\label{eq:SZ}
    |Z(F) \cap (A \times B \times C)| \ll_d n^2.
\end{equation}
This bound cannot be improved in general. 
For example, if $F(x,y)=x+y+z$, $A=B=\{1,\dots,n\}$,  and $C=\{-1,\dots,-n\}$,
then $|Z(F) \cap (A \times B \times C)|\gg n^2$. 
More generally,
if the equation $F(x,y,z)=0$ is in some sense equivalent to an equation of the form 
$\varphi_1(x)+\varphi_2(y)+\varphi_3(z) = 0$,
then we can choose $A,B,C$ so that $|Z(F) \cap (A \times B \times C)| \gg n^2$.
The following definition makes this property precise.

\begin{definition}\label{def:degenerate}
A polynomial $F(x,y,z)\in\R[x,y,z]$ is \emph{degenerate} if there are intervals $I_1,I_2,I_3$,
and for each $i$ there is a smooth (infinitely differentiable) function $\varphi_i:I_i\to \R$ which has a smooth inverse, 
such that for all $(x,y,z)\in I_1\times I_2\times I_3$ we have $F(x,y,z)=0$ if and only if $\varphi_1(x)+\varphi_2(y)+\varphi_3(z) = 0$.
\end{definition}

Elekes and Szab\'o \cite{ES12} showed that if the  polynomial is not degenerate in this sense, then the bound \eqref{eq:SZ} can be improved to $n^{2-\eta}$ for some $\eta>0$.
A quantitative improvement to $\eta = 1/6$ was obtained by Raz, Sharir and de Zeeuw \cite{RSZ16}, leading to the following statement.

\begin{theorem}[\cite{ES12,RSZ16}]\label{thm:ES}
Let $F\in \R[x,y,z]$ be a polynomial of degree $d$.
If $F$ is not degenerate, 
then for any $A,B,C\subset \R$ of size $n$ we have
\[
\left|Z(F)\cap (A\times B\times C)\right| \ll_d n^{2-1/6}.
\]
\end{theorem}

Not much attention has been paid to lower bound constructions for this theorem.
Elekes \cite{E99} noted that for $F = x^2+xy+y^2-z$ and $A = \{1,\ldots,n\}$ we have $|Z(F)\cap(A\times A\times A)| \gg n\sqrt{\log{n}}$ (actually, Elekes formulated this in a different way, which we mention in the next section; see \cite{Z18} for more discussion).
This was the only known lower bound for Theorem \ref{thm:ES},
and some have suggested that the upper bound could be improved as far as $O(n^{1+\varepsilon})$ for an arbitrarily small $\varepsilon >0$;
for instance, the fourth author wrote this in \cite{Z18}.
 
The main purpose of this paper is to show by means of a simple example that this is not the case, 
and that in fact the bound in Theorem \ref{thm:ES} cannot be improved beyond $O(n^{3/2})$.
Our main result is the following theorem.

\begin{theorem}\label{thm:main1} 
There exists a polynomial $F\in \R[x,y,z]$ of degree $2$ that is not degenerate, 
such that for any $n$ there is a set $A\subset \R$ of size $n$ with
\[
\left|Z(F)\cap (A \times A \times A)\right| \gg n^{3/2}.
\]
\end{theorem}

In Section \ref{sec:morevars},
we briefly discuss possible extensions of this theorem to polynomials in more variables.

\subsection{The Elekes-R\'{o}nyai Problem}

Before the work of Elekes and Szab\'o \cite{ES12}, 
Elekes and R\'onyai \cite{ER00} considered the question of bounding the image of a polynomial $f\in \R[x,y]$ restricted to a Cartesian product, assuming that $f$ does not have a certain special form, which is specified in the following definition.

\begin{definition}\label{def:addmult}
A polynomial $f(x,y)\in\R[x,y]$ is \emph{additive} if there are polynomials $g,h,k\in \R[t]$ such that $f(x,y) = g(h(x)+k(y))$,
and it is \emph{multiplicative} if there are polynomials $g,h,k\in \R[t]$ such that $f(x,y) = g(h(x)\cdot k(y)$.
\end{definition}

Elekes and R\'onyai \cite{ER00} proved that if $f\in \R[x,y]$ is not additive or multiplicative,
then for every $A,B \subseteq \R$ with $|A|=|B|=n$ the image $|f(A,B)|$ is superlinear in $n$.
The current state of the art for this problem is the following result of Raz, Sharir and Solymosi \cite{RSS16}.

\begin{theorem}[\cite{ER00, RSS16}]\label{thm:ER}
Let $f\in \R[x,y]$ be a polynomial of degree $d$.
If $f$ is not additive or multiplicative,
then for any $A,B\subset \R$ of size $n$ we have
\[|f(A,B)| \gg_d n^{4/3}.\]
\end{theorem}

Elekes \cite{E99} noted that if $f(x,y) = x^2+xy+y^2$ and $A = \{1,\ldots, n\}$,
then $|f(A,A)| \ll n^2/\sqrt{\log{n}}$.
This is the best known upper bound construction for Theorem \ref{thm:ER}, which suggests that we may have $|f(A,B)| \gg n^{2-\epsilon}$ for all positive $\epsilon$.
This conjecture is widely believed, see for instance Elekes \cite{E99} or Matou\v sek \cite[Section 4.1]{M02}.
The construction that we give in the proof of Theorem \ref{thm:main1} does not translate into a construction that disproves this conjecture.
 
Nevertheless,
we show that 
there is a polynomial that takes only a linear number of values on a certain large subset of the pairs in $A\times A$.
This approach is partly inspired by work of Alon, Ruzsa and Solymosi \cite{ARS18} concerning constructions for the sum-product problem along graphs. See also \cite{RNW18} for a slightly improved construction.

Let $G$ be a bipartite graph on $A$ and $B$ with edge set $E(G)\subset A\times B$.
For a polynomial $f\in \R[x,y]$ we define the image of $f$ along $G$ to be  $f_G(A,B)= \left\{f(a,b): (a,b) \in E(G) \right\}$.
Our result is the following.
\begin{theorem}\label{thm:main2}
There exists a polynomial $f\in\R[x,y]$ of degree $2$ that is not additive or multiplicative,
a finite set $A\subset \R$ of size $n$,
and a bipartite graph $G$ on $A\times A$, such that 
\[|E(G)|\gg n^{3/2}~~~~~~\text{and}~~~~~~
|f_G (A,B) | \leq n.\]
\end{theorem}

\section{The Elekes--Szab\'{o} problem}

In this section we prove Theorem \ref{thm:main1}.

Define
\[F(x,y,z) = (x-y)^2 + x - z.\] 
We set $A = \{1,\ldots,n\}$ and we consider the intersection of $F$ with $A\times A\times A$. 
Consider the subset
\[
T = \left\{(k,k+\ell,k+\ell^2): 
k,\ell\in \Z, ~
0 \leq k \leq n/2,  ~
0 \leq \ell \le \sqrt{n}/2 \right\} \subset A\times A\times A.
\]	
Each choice of $k$ and $\ell$ determines a distinct triple in $T$, and
so we have $|T| \gg n^{3/2}$.	
For each triple in $T$, we have 
\[F(k,k+\ell,k+\ell^2) = (k - (k+\ell))^2 + k -(k+\ell^2) = 0,\]
so $T\subset Z(F)$. 
Therefore we have
\[
\left|Z(F)\cap (A \times A \times A)  \right| \gg n^{3/2}.
\]

It remains to show that $F$ is not degenerate in the sense of Definition \ref{def:degenerate}.
We will use an idea introduced by Elekes and R\'onyai \cite{ER00},
which is that this type of degeneracy can be verified using the following straightforward derivative test; see for instance \cite[Lemma 33]{ES12} or \cite[Lemma 2.2]{Z18}.

\begin{lemma}\label{lem:difftest}
Let $f:\R^2\to \R$ be a smooth function on some open set $U\subset \R^2$ with $f_x$ and $f_y$ not identically zero.
If there exist smooth functions $\psi, \varphi_1,\varphi_2$ on $U$ such that 
\[ f(x,y) = \psi(\varphi_1(x) + \varphi_2(y)),\]
then
\begin{equation}\label{eq:test} \frac{\partial^2\left(\log|f_x/f_y|\right)}{\partial x\partial y}
\end{equation}
is identically zero on $U$.
\end{lemma}

Suppose that $F(x,y,z) = (x-y)^2+x-z$ is degenerate,
so in some neighborhood $I_1\times I_2\times I_3$ we have $F(x,y,z) = 0$ if and only if $\varphi_1(x)+\varphi_2(y) + \varphi_3(z) = 0$.
Then, since $\varphi_3$ has a smooth inverse on $I_3$, 
we can write $\psi(t)  = \varphi_3^{-1}(-t)$,
so that $F(x,y,z) = 0$ is equivalent to  $z = \psi(\varphi_1(x) + \varphi_2(y))$.
At the same time, $F(x,y,z) = 0$ rewrites to $z = (x-y)^2 + x$, so on $I_1\times I_2\times I_3$ we have
\[ \psi(\varphi_1(x) + \varphi_2(y)) = (x-y)^2 + x.\]

We now check if the expression \eqref{eq:test} is identically zero on $I_1\times I_2\times I_3$.
We have
\[\log|f_x/f_y| = \log\left|\frac{2(x-y)+1}{-2(x-y)}\right|
= \log|2x-2y+1| - \log|2x-2y|,\]
so
\[\frac{\partial^2\left(\log|f_x/f_y|\right)}{\partial x\partial y}
= \frac{\partial}{\partial x} \left(\frac{-1}{x-y+1/2} + \frac{1}{x-y}\right)
=\frac{1}{(x-y+1/2)^2} - \frac{1}{(x-y)^2}.
\]
This expression equals zero only when $y-x = 1/4$, so it does not vanish on any nontrivial open set.  
Thus \eqref{eq:test} is not identically zero, and by Lemma \ref{lem:difftest} this contradicts our assumption that $F$ is degenerate.

\section{The Elekes--R\'{o}nyai problem along a graph}

We now prove Theorem \ref{thm:main2}, concerning the image of a polynomial along a subset of a Cartesian product. 

Define the polynomial 
\[f(x,y) = (x - y)^2 + x.\] 
Set $A = \{1,\ldots,n\}$ and let $G$ be the bipartite graph on $A\times A$ with the edge set
\[E(G) = 
\left\{(k,k+\ell) :k,\ell\in \Z, ~
0 \leq k \leq n/2,  ~
0 \leq \ell \le \sqrt{n}/2
\right\}\subset A\times A.\]
We have $|E(G)| \gg n^{3/2}$.
Applying $f$ along any edge gives a non-negative integer
$$ (k-(k+\ell))^2 +k  < n.$$
This shows that
$$|f_G (A \times A)| \leq n.$$

It remains to prove that $f$ is not additive or multiplicative.
We could again do this using Lemma \ref{lem:difftest}, but here we can use a more elementary approach.
We treat the two cases separately.

\textit{Additive case:} Suppose $f(x,y) = g(h(x)+k(y))$.  
Note that $g$, $h$ and $k$ must have degree at most $2$.
We cannot have $\deg(g) = 1$,
since then $f(x,y)$ would not have any cross term $xy$.
If $\deg (g)=2$, then $\deg(h) = \deg(k) = 1$.
We can write
\[g(t)=a_2t^2+a_1t+a_0,
~~~ h(x)=b_1x+b_0,
~~~k(y)=c_1y+c_0,\]
with $b_1$ and $c_1$ non-zero.
Then we have
\begin{equation}\label{eq:long}
f(x,y)=(x-y)^2+x=a_2(b_1x+b_0+c_1y+c_0)^2+a_1(b_1x+b_0+c_1y+c_0)+a_0.
\end{equation}
Calculating the coefficient for the $y$ term on the right hand side and comparing with the left hand side, it follows that
\begin{equation}
\label{contradiction}
    2a_2(b_0+c_0)+a_1=0.
\end{equation}
On the other hand, calculating the coefficient for the $x$ term on the right hand side of \eqref{eq:long} and comparing with the left hand side, it follows that 
\[b_1(2a_2(b_0+c_0)+a_1)=1.\]
Since $b_1\neq 0$, this contradicts \eqref{contradiction}. 

\textit{Multiplicative case:} Suppose $f(x,y) = g(h(x)\cdot k(y))$. 
We cannot have  $\deg(g)=2$, since then $h$ or $k$ would have to be constant, and $f(x,y)$ would not depend on both variables. 
Therefore we have $\deg(g)=1$. 
In this case, we must have $\deg(h)=\deg(k) = 1$.
We can write
\[g(t)=a_1t+a_0,~~~h(x)=b_1x+b_0,~~~  k(y)=c_1y+c_0\]
and
\[
(x-y)^2+x=f(x,y)=a_1((b_1x+b_0)(c_1y+c_0))+a_0.
\]
This is a contradiction, since there is no $x^2$ or $y^2$ term on the right hand side.

This completes our proof that $f$ is not additive or multiplicative,
which completes our proof of Theorem \ref{thm:ER}.

\section{Extensions to more variables}\label{sec:morevars}

\subsection{Four variables}

One can consider the same problems for polynomials in more variables.
Raz, Sharir and de Zeeuw \cite{RSZ18} proved that
for $F\in \R[x,y,s,t]$ of degree $d$ and $A,B,C,D\subset \R$ of size $n$,
we have
\begin{equation}\label{eq:ES4D}
\left|Z(F)\cap (A\times B\times C\times D)\right| \ll_d n^{8/3},
\end{equation}
unless $F(x,y,s,t)=0$ is in a local sense (similar to Definition \ref{def:degenerate}) equivalent to an equation of the form $\varphi_1(x)+\varphi_2(y)+\varphi_3(s)+\varphi_4(t) = 0$.

A construction of Valtr \cite{V05} (see also \cite[Section 5.3]{SSZ13}) essentially shows that for
\[V(x,y,s,t) = (x-y)^2+s-t\]
one can set $A = B = \{1,\ldots,n^{2/3}\}$ and $C = D = \{1,\ldots, n^{4/3}\}$, so that
\[
\left|Z(V)\cap (A\times B\times C\times D)\right| \gg n^{8/3}.
\]
This would show that \eqref{eq:ES4D} is tight, if it weren't for the fact that $A,B$ and $C,D$ have different sizes.
(A similar, older, construction of Elekes \cite[Example 1.16]{E02} achieves the same with the polynomial $xy+s-t$, but is less relevant to us here.)

If we require that $A,B,C,D$ have the same size (and then we may as well assume that they all equal $A\cup B\cup C\cup D$),
then we can take Valtr's polynomial $V(x,y,s,t)$ together with the set $A =\{1,\ldots, n\}$. 
Similarly to in our proof of Theorem \ref{thm:main1},
considering quadruples of the form $(k,k+\ell, m, m+\ell^2)$ with $1\leq k,m\leq n/2$ and $ 1\leq \ell\leq \sqrt{n}/2$,
we get 
\[|Z(V)\cap(A\times A\times A\times A)| \gg n^{5/2}.\]
It is not hard to verify (as in our proof of Theorem \ref{thm:main1}) that $V(x,y,s,t)$ is not degenerate in the sense of \cite{RSZ18},
so this gives a lower bound construction for \eqref{eq:ES4D},
which is the best known.

Note that the polynomial $F(x,y,z)$ in our proof of Theorem \ref{thm:main1} can be obtained from Valtr's polynomial $V(x,y,s,t)$ by setting $s=x$ and $t=z$.

\subsection{More than four variables}

For more than four variables, we do not have a statement that is entirely analogous to Theorem \ref{thm:ES} or \eqref{eq:ES4D}.
Bays and Breuillard \cite{BB18} proved a similar statement for any number of variables, but without an explicit exponent, and with a different description of the exceptional form.
Also, Raz and Tov \cite{RT18} extended Theorem \ref{thm:ER} to any number of variables, with an explicit exponent.

Because for the Elekes--Szab\'o problem in more than four variables we do not have explicit exponents, and also because the appropriate definition of degeneracy is not clear,
we only briefly touch on constructions for more variables here.

There are various ways of extending our constructions to more variables;
 one can for instance take the polynomial 
\[F(x_1,\dots,x_m)= (x_1+\dots+x_{m-1})^2+x_1-x_m
\]
and the grid $A^m$, where $A = [-n,2n]$. Consider the set
$$T = \left\{ 
(k_1, k_2 - k_1, 
\dots , k_{m-2}-k_{m-3}, \ell - k_{m-2}, k_1 +
\ell^2) : 0 \leq k_i \leq n, \ 0 \leq \ell \leq \sqrt{n}  \right\}.$$
Then we have $T\subset Z(F)\cap A^m$, 
which implies
\[
\big|Z(F)\cap A^m  \big| \gg n^{m - \frac{3}{2}}.
\]
This should be compared with the Schwartz--Zippel bound $|Z(F)\cap A^m|\ll n^{m-1}$.
A potential Elekes--Szab\'o theorem in $m$ variables, i.e. an explicit version of the result of Bays and Breuillard,
would give a bound of the form $|Z(F)\cap A^m| \ll n^{m-1-\eta_m}$ for some $\eta_m>0$, 
under the condition that $F$ is not degenerate in some sense.
Presuming that our polynomial $F$ is not of this form, it would show that we must have $\eta_m\leq 1/2$.

\section*{Acknowledgements} 

Mehdi Makhul was supported by the Austrian Science Fund (FWF): W1214-N15, Project DK9. Oliver Roche-Newton and Audie Warren were partially supported by the Austrian Science Fund FWF Project P 30405-N32. We are very grateful to Niels Lubbes, for a particularly interesting conversation which resulted in us attempting to find non-trivial constructions for the Elekes--Szab\'{o} problem. Thanks also to J\'ozsef Solymosi for pointing out some helpful and relevant references.


\end{document}